# NONPARAMETRIC ESTIMATION OF MEAN-SQUARED PREDICTION ERROR IN NESTED-ERROR REGRESSION MODELS

By Peter Hall and Tapabrata Maiti[1]

*Australian National University and Iowa State University*

Nested-error regression models are widely used for analyzing clustered data. For example, they are often applied to two-stage sample surveys, and in biology and econometrics. Prediction is usually the main goal of such analyses, and mean-squared prediction error is the main way in which prediction performance is measured. In this paper we suggest a new approach to estimating mean-squared prediction error. We introduce a matched-moment, double-bootstrap algorithm, enabling the notorious underestimation of the naive mean-squared error estimator to be substantially reduced. Our approach does not require specific assumptions about the distributions of errors. Additionally, it is simple and easy to apply. This is achieved through using Monte Carlo simulation to implicitly develop formulae which, in a more conventional approach, would be derived laboriously by mathematical arguments.

**1. Introduction.** Unbalanced nested-error regression models often arise in two-stage sample surveys, multilevel modeling, biological experiments and econometric analysis. Beside the noise, a source of variation is added to explain the correlation among observations within clusters, or subjects, and to allow the analysis to borrow strength from other clusters. Such nested-error regression models are particular cases of general linear mixed models, which often form the basis for inference about small-area means or subject-specific values.

In this article we propose a new, nonparametric bootstrap technique for estimating the mean-squared error of predictors of mixed effects. The new

Received March 2005; revised September 2005.
[1]Supported in part by NSF Grant SES-03-18184.

method has several attractive properties. First, it does not require specific distributional assumptions about error distributions. Second, it produces positive, bias-corrected estimators of mean-squared prediction errors. (See [2] and [5] for discussion of possible negativity.) Third, it is easy to apply. Although our emphasis is on small-area prediction, our methodology is equally useful for other applications, such as estimating subject- or cluster-specific random effects.

Standard mixed-effects prediction involves two steps. First, a best linear unbiased predictor, or BLUP, is derived under the assumption that model parameters are known. Then, the model parameters are replaced by estimators, producing an empirical version of BLUP. This approach is popular because it is straightforward and, at this level, does not require distributional assumptions.

However, estimation of mean-squared prediction error is significantly more challenging. The variability of parameter estimators can substantially influence mean-squared error, to a much greater extent than a conventional asymptotic analysis suggests. Moreover, the nature and extent of this influence is intimately connected to the values of the design variables and to properties of the two error distributions.

In this paper we point out that, in terms of the biases of estimators of mean-squared prediction error, the two error distributions influence results predominantly through their second and fourth moments. This observation leads to a surprisingly simple, moment-matching, double-bootstrap algorithm for estimating, and correcting for, bias. We show that this approach substantially reduces the large degree of underestimation by the naive approach.

Kackar and Harville [20] and Harville and Jeske [18] studied various approximations to the mean-squared prediction error of the empirical BLUP, assuming normality in both stages. Prasad and Rao [27] pointed out that if unknown model parameters are replaced by their estimators, then significant underestimation of true mean-squared prediction error can still result. This difficulty can have significant impact on policy making. To alleviate it, Prasad and Rao [27] constructed second-order correct mean-squared error estimators under normal models. Datta and Lahiri [8] extended the Prasad–Rao approach to cases where model parameters are estimated using maximum likelihood, or restricted maximum likelihood, methods. Das, Jiang and Rao [6] gave rigorous proofs of these results under normality. Bootstrap methods in parametric settings have been suggested, for this problem, by Booth and Hobert [3] and Lahiri [22], for example.

Jiang, Lahiri and Wan [19] proposed a jackknife-based bias correction of the mean-squared error estimator. Again, unlike the approach taken in the present paper, explicit parametric models are required. The problem of mixed-effects prediction is in part one of deconvolution, and so conventional, nonparametric jackknife estimators of mean-squared error are not applicable; hence the need by Jiang, Lahiri and Wan [19] for parametric assumptions. For convenient implementation the methods proposed there also require a closed-form expression to be available for the leading term in an expansion of mean-squared prediction error, as a function of unknown parameters. The main advantage of our technique is that it does not require parametric assumptions about the distributions of the two sources of error in the model, or an analytical partition of those sources.

On the other hand, the jackknife approach has advantages. Principal among these are the fact that it can be used beyond the linear setting treated in the present paper, for example in the case of generalized linear mixed models; and that, in a parametric context, related methods might potentially be employed to construct prediction intervals, rather than estimators of mean-squared prediction error. In the first of these settings, our method is confounded by nonlinear link functions. In the second, aspects of the error distributions, beyond just low-order moments, play necessary roles in constructing the prediction interval, and so again our moment-matching bootstrap approach is not suitable. Further discussion of jackknife methods in the small-area estimation problem is given by Lahiri [23].

An approach alternative to that given in this paper would be to estimate the two error distributions explicitly, and base a bootstrap algorithm on those estimators. However, since we wish to treat both error distributions nonparametrically, then the deconvolution problem would be quite nonstandard; the large literature on nonparametric deconvolution is devoted almost entirely to the case where one distribution is assumed known and the other is estimated. Early work in this area includes that of Carroll and Hall [4] and Fan [12, 13], and more recent contributions can be accessed through citations by, for example, Delaigle and Gijbels [9].

Identifiability of the full, doubly nonparametric deconvolution problem rests on the fact that some of the measurements are repeated. Methods for solution can be developed, for example, by starting from the approaches introduced by El-Amraoui and Goffinet [11] and Li and Vuong [26] in different contexts. However, in addition to the intrinsic difficulty, to a practitioner, of implementing a full deconvolution approach, such a technique would involve choosing smoothing parameters, which would have to be selected to optimize performance in a nonstandard problem where the target is bias reduction, not density estimation. By way of comparison, the bootstrap approach suggested in the present paper is simple and explicit. Only low-order moment estimators of the error distributions are required, and the estimators are directly defined as functions of the data.

## 2. Methodology.

2.1. *Model.* We observe data pairs $(X_{ij}, Y_{ij})$ generated by the model

$$(2.1) \quad Y_{ij} = \mu + X'_{ij}\beta + U_i + s_{ij}V_{ij} \qquad \text{for } 1 \le i \le n \text{ and } 1 \le j \le n_i,$$

where each $n_i \ge 2$, $Y_{ij}$ and $\mu$ are scalars, $X_{ij}$ is an $r$-vector, $\beta$ is an $r$-vector of unknown parameters, the scalar $s_{ij}$ is known (generally as a function of $X_{i1}, \dots, X_{in_i}$), the $U_i$'s and $V_{ij}$'s are totally independent, the $U_i$'s are identically distributed, the $V_{ij}$'s are identically distributed, $E(U_i) = E(V_{ij}) = 0$ for each $i, j$, $E(U_i^2) = \sigma_U^2$ and $E(V_{ij}^2) = \sigma_V^2$. All inference will be conducted conditionally on $\mathcal{X}$, which denotes the set of explanatory data $X_{ij}$ for $1 \le i \le n$ and $1 \le j \le n_i$.

The model (2.1) is a generalization of the unbalanced nested-error regression model [29, 31], and is commonly used to model two-level clustered data. For example, Battese, Harter and Fuller [1] and Datta and Ghosh [7] used this model, with $s_{ij} \equiv 1$, for predicting the areas under corn and soybeans for 12 counties in North Central Iowa. Rao and Choudhry [30] studied the population of unincorporated tax filers from the province of Nova Scotia, Canada using (2.1) with $s_{ij} = X_{ij}^{1/2}$.

Of course, (2.1) arises through noise, in terms of the $V_{ij}$'s, being added to an observation,

$$\Theta_i = \mu + \underline{X}'_i\beta + U_i,$$

of the small-area modeling "parameter." Here $\underline{X}_i = n_i^{-1}\sum_j X_{ij}$. Our objective is to make inference about estimators of the performance of predictors of the small-area mean $\Theta_i$, or even just the random effect $U_i$ (in the case $\mu = 0$ and $\beta = 0$).

2.2. *Formulae for predictors.* Put $\bar{X}_i = a_i^{-1}\sum_j s_{ij}^{-2}X_{ij}$ and $\bar{Y}_i = a_i^{-1} \times \sum_j s_{ij}^{-2}Y_{ij}$, where $a_i = \sum_j s_{ij}^{-2}$. The best linear unbiased predictor of $\Theta_i$ is

$$\Theta_i^{\text{BLUP}} = \mu + \underline{X}'_i\beta + \rho_i(\bar{Y}_i - \mu - \bar{X}'_i\beta),$$

where $\rho_i = \sigma_U^2/(\sigma_U^2 + a_i^{-1}\sigma_V^2)$. Replacing $\mu$ and $\beta$ by their weighted least-squares estimators, $\tilde{\mu}$ and $\tilde{\beta}$ say, defined under the temporary assumption that $\sigma_U^2$ and $\sigma_V^2$ are known, we obtain an empirical version of $\Theta_i^{\text{BLUP}}$,

$$\widetilde{\Theta}_i^{\text{BLUP}} = \tilde{\mu} + \underline{X}'_i\tilde{\beta} + \rho_i(\bar{Y}_i - \tilde{\mu} - \bar{X}'_i\tilde{\beta}).$$

Here,

$$\tilde{\mu} = \left(\sum_{i=1}^n \mathbf{1}'_i\mathbf{W}_i^{-1}\mathbf{1}_i\right)^{-1}\sum_{i=1}^n \mathbf{1}'_i\mathbf{W}_i^{-1}(Y_i - \mathbf{X}'_i\tilde{\beta}),$$

$$\tilde{\beta} = \left\{\sum_{i=1}^n (\mathbf{X}'_i - \mathbf{1}_i\bar{X}')'\mathbf{W}_i^{-1}(\mathbf{X}'_i - \mathbf{1}_i\bar{X}')\right\}^{-1}\sum_{i=1}^n (\mathbf{X}'_i - \mathbf{1}_i\bar{X}')'\mathbf{W}_i^{-1}(Y_i - \bar{Y}\mathbf{1}_i),$$

where $\mathbf{1}_i$ is the vector of 1's of length $n_i$, $\mathbf{X}_i$ denotes the $r \times n_i$ matrix with $X_{ij}$ as its $j$th column, $\mathbf{W}_i$ is the $n_i \times n_i$ matrix of which the $(j_1, j_2)$th component is $\sigma_U^2 + \delta_{j_1 j_2} s_{ij_1}^2 \sigma_V^2$, $\delta_{j_1 j_2}$ is the Kronecker delta, $Y_i$ is the $n_i$-vector with $j$th component $Y_{ij}$, and

$$\bar{X} = \left( \sum_{i=1}^n \mathbf{1}_i' \mathbf{W}_i^{-1} \mathbf{1}_i \right)^{-1} \sum_{i=1}^n \mathbf{X}_i \mathbf{W}_i^{-1} \mathbf{1}_i,$$

$$\bar{Y} = \left( \sum_{i=1}^n \mathbf{1}_i' \mathbf{W}_i^{-1} \mathbf{1}_i \right)^{-1} \sum_{i=1}^n Y_i' \mathbf{W}_i^{-1} \mathbf{1}_i,$$

denote an $r$-vector and a scalar, respectively.

A practical form of $\widehat{\Theta}_i^{\mathrm{BLUP}}$ is

$$\widehat{\Theta}_i^{\mathrm{BLUP}} = \hat{\mu} + \underline{X}_i' \hat{\beta} + \hat{\rho}_i (\bar{Y}_i - \hat{\mu} - \bar{X}_i' \hat{\beta}),$$

where $\hat{\mu}$, $\hat{\beta}$ and $\hat{\rho}_i$ differ from $\tilde{\mu}$, $\tilde{\beta}$ and $\rho_i$, respectively, in that $\sigma_U^2$ and $\sigma_V^2$ are replaced by estimators, $\hat{\sigma}_U^2$ and $\hat{\sigma}_V^2$ say. We wish to construct a bias-corrected estimator of the mean-squared prediction error,

$$(2.2) \qquad \mathrm{MSE}_i = E\{ (\widehat{\Theta}_i^{\mathrm{BLUP}} - \Theta_i)^2 \mid \mathcal{X} \}.$$

2.3. *Formulae for $\hat{\sigma}_U^2$ and $\hat{\sigma}_V^2$.* The estimators used here are borrowed from [31]. Put

$$(2.3) \qquad \bar{V}_i = \frac{\sum_j s_{ij}^{-1} V_{ij}}{\sum_j s_{ij}^{-2}}, \qquad \bar{X}_i = \frac{\sum_j s_{ij}^{-2} X_{ij}}{\sum_j s_{ij}^{-2}}, \qquad \bar{Y}_i = \frac{\sum_j s_{ij}^{-2} Y_{ij}}{\sum_j s_{ij}^{-2}},$$

$p_{ij} = s_{ij}^{-1} (X_{ij} - \bar{X}_i)$, $q_{ij} = s_{ij}^{-1} (Y_{ij} - \bar{Y}_i)$ and $e_{ij} = V_{ij} - \bar{V}_i$. In this notation,

$$(2.4) \qquad q_{ij} = p_{ij}' \beta + e_{ij}, \qquad 1 \le j \le n_i, 1 \le i \le n.$$

Note too that $E(e_{ij}) = 0$ and $\mathrm{cov}(e_{ij_1}, e_{ij_2}) = t_{ij_1 j_2} \sigma_V^2$, where

$$t_{ij_1 j_2} = \delta_{j_1 j_2} - \frac{s_{ij_1}^{-1} + s_{ij_2}^{-1} - 1}{\sum_j s_{ij}^{-2}}.$$

Let $\mathbf{P}_i = (p_{i1}, \ldots, p_{i,n_i-1})$ be an $r \times (n_i-1)$ matrix, and let $\mathbf{P} = (\mathbf{P}_1, \ldots, \mathbf{P}_n)$ be an $r \times (N-n)$ matrix, where $N = \sum_i n_i$. Let $q_i = (q_{i1}, \ldots, q_{i,n_i-1})'$ and $e_i = (e_{i1}, \ldots, e_{i,n_i-1})'$ be $(n_i-1)$-vectors, and let $q = (q_1', \ldots, q_n')'$ and $e = (e_1', \ldots, e_n')'$ be $(N-n)$-vectors. Let $\mathbf{T}_i$ be the $(n_i-1) \times (n_i-1)$ matrix with $t_{ij_1 j_2}$ in position $(j_1, j_2)$, and let $\mathbf{T}$ be the $(N-n) \times (N-n)$ matrix with blocks $\mathbf{T}_1, \ldots, \mathbf{T}_n$ down the main diagonal and blocks of zeros elsewhere. Bearing in mind linear relationships, the set of equations (2.4) is equivalent to

$$(2.5) \qquad q = \mathbf{P}' \beta + e,$$

where $E(e) = 0$ and $\mathrm{cov}(e) = \mathbf{T}\sigma_V^2$.

We shall assume that $N - n > r$, and that the matrices $\mathbf{T}$ and $\mathbf{P}\mathbf{T}^{-1}\mathbf{P}'$ are both of full rank, $N - n$ and $r$, respectively. Then, the sum of squares for error arising from the regression model (2.5) is $\mathrm{SSE}_1 = \hat{e}'\mathbf{T}^{-1}\hat{e}$, where $\hat{e} = q - \mathbf{P}'\hat{\beta}^{\mathrm{wls}}$ is the vector of residuals and $\hat{\beta}^{\mathrm{wls}} = (\mathbf{P}\mathbf{T}^{-1}\mathbf{P}')^{-1}\mathbf{P}\mathbf{T}^{-1}q$ is the weighted least-squares estimator of $\beta$. It is shown in a longer version of this paper [17] that, under the full-rank conditions,

$$(2.6) \qquad E(\mathrm{SSE}_1) = (N - n - r)\sigma_V^2.$$

This property motivates the estimator

$$(2.7) \qquad \hat{\sigma}_V^2 = \frac{\mathrm{SSE}_1}{N - n - r}.$$

In (2.6) and below we interpret expected value to be taken conditional on the set $\mathcal{X}$ of design variables. That is, we drop the notation "$|\mathcal{X}$" used at (2.2).

Write $\bar{p}_{ij} = s_{ij}^{-1}X_{ij}$, $\bar{q}_{ij} = s_{ij}^{-1}Y_{ij}$ and $\bar{e}_{ij} = s_{ij}^{-1}U_i + V_{ij}$ for the uncentred versions of $p_{ij}$, $q_{ij}$ and $e_{ij}$, and put $\bar{t}_{ij_1 j_2} = (s_{ij_1}s_{ij_2})^{-1}\sigma_U^2 + \delta_{j_1 j_2}\sigma_V^2$. Let $\bar{\mathbf{P}}_i = (p_{i1}, \ldots, p_{i,n_i})$ be an $r \times n_i$ matrix, let $\bar{\mathbf{P}} = (\bar{\mathbf{P}}_1, \ldots, \bar{\mathbf{P}}_n)$ be an $r \times N$ matrix, let $\bar{q}_i = (\bar{q}_{i1}, \ldots, \bar{q}_{in_i})'$ and $\bar{e}_i = (\bar{e}_{i1}, \ldots, \bar{e}_{i,n_i})'$ be $n_i$-vectors, and let $\bar{q} = (\bar{q}_1', \ldots, \bar{q}_n')'$ and $\bar{e} = (\bar{e}_1', \ldots, \bar{e}_n')'$ be $N$-vectors. Let $\bar{\mathbf{T}}_i$ be the $n_i \times n_i$ matrix with $\bar{t}_{ij_1 j_2}$ in position $(j_1, j_2)$, and let $\bar{\mathbf{T}}$ be the $N \times N$ block-diagonal matrix with blocks $\bar{\mathbf{T}}_1, \ldots, \bar{\mathbf{T}}_n$ down the main diagonal. In this notation, the model at (2.1) is equivalent to

$$(2.8) \qquad \bar{q} = \bar{\mathbf{P}}'\beta + \bar{e},$$

where $E(\bar{e}) = 0$ and $\mathrm{cov}(\bar{e}) = \bar{\mathbf{T}}$.

Assuming $\bar{\mathbf{P}}$ is of full rank, $r$, the sum of squares for error arising from (2.8) is $\mathrm{SSE}_2 = \hat{\bar{e}}'\hat{\bar{e}}$, where $\hat{\bar{e}} = \bar{q} - \bar{\mathbf{P}}'\hat{\beta}^{\mathrm{ols}}$ is a new vector of residuals, and $\hat{\beta}^{\mathrm{ols}} = (\bar{\mathbf{P}}\bar{\mathbf{P}}')^{-1}\bar{\mathbf{P}}\bar{q}$ is the ordinary least-squares estimator. In a longer version of this paper [17] it is proved that, analogously to (2.6),

$$(2.9) \qquad E(\mathrm{SSE}_2) = K\sigma_U^2 + (N - r)\sigma_V^2,$$

where $K = K_1 - K_2$, $K_1 = \sum_i \sum_j s_{ij}^{-2}$ and

$$K_2 = \sum_{i=1}^n \left( \sum_{j=1}^{n_i} s_{ij}^{-2}X_{ij} \right)^{\mathrm{T}} \left( \sum_{i=1}^n \sum_{j=1}^{n_i} s_{ij}^{-2}X_{ij}X_{ij}' \right)^{-1} \sum_{j=1}^{n_i} s_{ij}^{-2}X_{ij}.$$

Property (2.9) suggests the estimator

$$(2.10) \qquad \hat{\sigma}_U^2 = \max[K^{-1}\{\mathrm{SSE}_2 - (N - r)\hat{\sigma}_V^2\}, 0].$$

Recall that the estimators $\hat{\sigma}_U^2$ and $\hat{\sigma}_V^2$ are substituted for $\sigma_U^2$ and $\sigma_V^2$, respectively, in formulae for $\tilde{\mu}$ and $\tilde{\beta}$, to obtain the estimators $\hat{\mu}$ and $\hat{\beta}$, respectively. In particular, they are substituted for $\sigma_U^2$ and $\sigma_V^2$ in the formula $\mathbf{W}_i = \sigma_U^2 I_{n_i} + \sigma_V^2 \operatorname{diag}(s_{i1}^2, \ldots, s_{in_i}^2)$, to obtain $\widehat{\mathbf{W}}_i$ say, and we need to invert $\widehat{\mathbf{W}}_i$ when computing $\hat{\mu}$ and $\hat{\beta}$. In some problems, a realistic discrete model for $U$ and $V$ can involve both taking the value zero with nonzero probability, and in this case there is a nonzero probability that $\widehat{\mathbf{W}}_i^{-1}$ is not well defined. More generally, there might be concern about cases where the determinant of $\widehat{\mathbf{W}}_i$ is positive but close to zero. To remove these theoretical pathologies it is sufficient to replace $\mathrm{SSE}_1$ by $\max(\mathrm{SSE}_1, \delta_n)$, where $\delta_n > 0$ denotes a small ridge parameter. See Section 4 for further discussion.

2.4. *Expansions, and analytical estimators, of* $\mathrm{MSE}_i$. Recall the definition of $\mathrm{MSE}_i$ at (2.2). It can be proved that $\hat{\Theta}_i^{\mathrm{BLUP}} - \Theta_i = \Delta_i + O_p(n^{-1/2})$, where $\Delta_i = \rho_i \bar{V}_i - (1 - \rho_i) U_i$ and $\bar{V}_i$ is as at (2.3). This property suggests that $\mathrm{MSE}_i = E(\Delta_i^2) + O(n^{-1})$. Indeed,

$$(2.11) \qquad \mathrm{MSE}_i = E(\Delta_i^2) + n^{-1}\psi_1(\xi_1) + O(n^{-2}),$$

where, here and below, $\psi_j$ denotes a smooth function depending only on the known design variables $X_{ij}$ and standard deviations $s_{ij}$, and $\xi_1 = (\sigma_U^2, \sigma_V^2, EU^4, EV^4)$ is the vector consisting of second and fourth moments of $U$ and $V$. See (4.2) in Theorem 1 in Section 4.1 for a rigorous formulation of (2.11).

It is readily seen that

$$(2.12) \qquad E(\Delta_i^2) = \frac{\sigma_U^2 a_i^{-1} \sigma_V^2}{\sigma_U^2 + a_i^{-1}\sigma_V^2},$$

where $a_i = \sum_j s_{ij}^{-2}$ is as in Section 2.2. Therefore, (2.11) can be written as

$$(2.13) \qquad \mathrm{MSE}_i = \psi_0(\xi_0) + n^{-1}\psi_1(\xi_1) + O(n^{-2}),$$

where $\xi_0 = (\sigma_U^2, \sigma_V^2)$ and $\psi_0$ is another known, smooth function.

From (2.13) it can be appreciated that, in order to estimate $\mathrm{MSE}_i$, we need only compute estimators of the second and fourth moments of $U$ and $V$, and substitute them into the approximate formula, $\mathrm{MSE}_i \approx \psi_0(\xi_0) + n^{-1}\psi_1(\xi_1)$. However, this will introduce a bias of size $n^{-1}$, since if $\hat{\xi}_0 = (\hat{\sigma}_U^2, \hat{\sigma}_V^2)$, then

$$(2.14) \qquad E\{\psi_0(\hat{\xi}_0)\} = \psi_0(\xi_0) + n^{-1}\psi_2(\xi_1) + O(n^{-2}),$$

where $\psi_2$ is a further known, smooth function. A rigorous formulation of (2.14) is given in (4.3) in Theorem 1.

In fact, if we take $\hat{\xi}_1$ to be a vector of root-$n$ consistent estimators of the respective components of $\xi_1$, then, since

$$(2.15) \qquad E\{\psi_1(\hat{\xi}_1)\} = \psi_1(\xi_1) + O(n^{-1}),$$

the estimator

$$(2.16) \qquad \widetilde{\mathrm{MSE}}_i = \psi_0(\hat{\xi}_0) + n^{-1}\psi_1(\hat{\xi}_1)$$

will satisfy

$$E(\widetilde{\mathrm{MSE}}_i) = \mathrm{MSE}_i + n^{-1}\psi_2(\xi_1) + O(n^{-2}). \tag{2.17}$$

We can correct for the term $n^{-1}\psi_2(\xi_1)$ on the right-hand side of (2.17) by moving it to the left, and replacing $\xi_1$ by its estimator,

$$E\{\widetilde{\mathrm{MSE}}_i - n^{-1}\psi_2(\hat{\xi}_1)\} = \mathrm{MSE}_i + O(n^{-2}). \tag{2.18}$$

Here we have used the fact that

$$E\{\psi_2(\hat{\xi}_1)\} = \psi_2(\xi_1) + O(n^{-1}). \tag{2.19}$$

[Result (4.4) in Theorem 1 gives (2.15) and (2.19) under explicit regularity conditions.] Property (2.18) suggests a bias-corrected estimator,

$$\widetilde{\mathrm{MSE}}_i^{\mathrm{bc}} = \widetilde{\mathrm{MSE}}_i - n^{-1}\psi_2(\hat{\xi}_1), \tag{2.20}$$

of $\mathrm{MSE}_i$, and argues that it has bias of order $n^{-2}$:

$$E(\widetilde{\mathrm{MSE}}_i^{\mathrm{bc}}) = \mathrm{MSE}_i + O(n^{-2}). \tag{2.21}$$

See Section 4.1 for discussion. Of course, (2.20) is motivated by the fact that the quantity

$$\widetilde{\mathrm{bias}}_i = n^{-1}\psi_2(\hat{\xi}_1) \tag{2.22}$$

is an estimator of the bias of $\widetilde{\mathrm{MSE}}_i$.

While the estimators at (2.16) and (2.20) might be satisfactory from a theoretical viewpoint, they are impractical or unattractive on several grounds. First, although the functions $\psi_1$ and $\psi_2$ are in principle known, they are very complicated functions of the $X_{ij}$'s and $s_{ij}$'s, and so implementing the estimators is not attractive to a practitioner. Second, the additive and subtractive nature of the corrections implicit in the procedures carries a risk that, in small to moderate samples, the estimators of $\mathrm{MSE}_i$ will be negative. Third, the complexity of the functions $\psi_1$ and $\psi_2$ would lead one to suspect that the procedures will be highly asymptotic in character. In particular, $n$ will have to be quite large before reasonably unbiased estimators will be obtained. Taken together, these difficulties motivate development of an alternative, bootstrap approach, which is likely to be more attractive. The bootstrap algorithm suggested below uses Monte Carlo simulation to approximate the functions $\psi_1$ and $\psi_2$, avoiding the need for explicit calculation.

2.5. *Bootstrap estimators of* $\mathrm{MSE}_i$. Results (2.13) and (2.14) imply that, in a bootstrap approach to this problem, it is sufficient from some viewpoints to resample from empirical "approximations" to the distributions of $U$ and $V$ that have first, second and fourth moments which are root-$n$ consistent for the corresponding moments of $U$ and $V$. In particular, we do not need the distributions from which we resample to actually be consistent for the distributions of $U$ and $V$. This is a variant of the moment-matching, or "wild," bootstrap method, which almost invariably addresses first, second and third, rather than first, second and fourth, moments. For recent applications of the moment-matching bootstrap, see [10, 14, 15, 16, 21, 25, 28].

With this motivation, we consider the following bootstrap algorithm. Given $z_2, z_4 > 0$ with $z_2^2 \leq z_4$, let $D(z_2, z_4)$ denote the distribution of a random variable $Z$, say, for which $E(Z) = 0$ and $E(Z^j) = z_j$ for $j = 2, 4$. Let $\mathcal{D}$ denote a class of such distributions, with exactly one member $D(z_2, z_4)$ for each pair $(z_2, z_4)$. Given the estimators $\hat{\sigma}_U^2$ and $\hat{\sigma}_V^2$ at (2.10) and (2.7), as well as estimators $\hat{\gamma}_U$ and $\hat{\gamma}_V$ of $\gamma_U = E(U^4)$ and $\gamma_V = E(V^4)$, satisfying the standard moment conditions $\hat{\sigma}_U^4 \leq \hat{\gamma}_U$ and $\hat{\sigma}_V^4 \leq \hat{\gamma}_V$, draw resamples $\mathcal{U}^* = \{U_1^*, \ldots, U_n^*\}$ and $\mathcal{V}^* = \{V_{ij}^* : 1 \leq i \leq n, 1 \leq j \leq n_i\}$ by sampling independently from the distributions $D(\hat{\sigma}_U^2, \hat{\gamma}_U)$ and $D(\hat{\sigma}_V^2, \hat{\gamma}_V)$, respectively, the distributions being the uniquely determined members of $\mathcal{D}$. Mimicking the model (2.1), define

$$Y_{ij}^* = \hat{\mu} + X_{ij}'\hat{\beta} + U_i^* + s_{ij}V_{ij}^* \qquad \text{for } 1 \leq i \leq n \text{ and } 1 \leq j \leq n_i.$$

Let $\mathcal{Z}$ and $\mathcal{Z}^*$ denote the set of all pairs $(X_{ij}, Y_{ij})$, and the set of all pairs $(X_{ij}, Y_{ij}^*)$, respectively. Using the data in $\mathcal{Z}^*$, compute the bootstrap versions $\hat{\mu}^*$, $\hat{\beta}^*$, $\hat{\sigma}_U^*$, $\hat{\sigma}_V^*$, $\hat{\gamma}_U^*$, $\hat{\gamma}_V^*$ and $\widehat{\Theta}_i^{*\text{BLUP}}$ of $\hat{\mu}$, $\hat{\beta}$, $\hat{\sigma}_U$, $\hat{\sigma}_V$, $\hat{\gamma}_U$, $\hat{\gamma}_V$ and $\widehat{\Theta}_i^{\text{BLUP}}$, respectively, and put

$$(2.23) \qquad \widehat{\text{MSE}}_i = E\{(\widehat{\Theta}_i^{*\text{BLUP}} - \Theta_i^*)^2 \mid \mathcal{Z}\};$$

compare (2.2). In (2.23), $\Theta_i^* = \hat{\mu} + \underline{X}_i'\hat{\beta} + U_i^*$. The quantity $\widehat{\text{MSE}}_i$ is our basic estimator of $\text{MSE}_i$. We shall prove in Section 4 that it has bias of order $n^{-1}$.

To bias-correct $\widehat{\text{MSE}}_i$ we use the double bootstrap, as follows. Conditional on $\mathcal{U}^*$ and $\mathcal{V}^*$, draw resamples $\{U_1^{**}, \ldots, U_n^{**}\}$ and $\{V_{ij}^{**} : 1 \leq i \leq n, 1 \leq j \leq n_i\}$ by sampling independently from the distributions $D\{(\hat{\sigma}_U^*)^2, \hat{\gamma}_U^*\}$ and $D\{(\hat{\sigma}_V^*)^2, \hat{\gamma}_V^*\}$, respectively. Let

$$Y_{ij}^{**} = \hat{\mu}^* + X_{ij}'\hat{\beta}^* + U_i^{**} + s_{ij}V_{ij}^{**} \qquad \text{for } 1 \leq i \leq n \text{ and } 1 \leq j \leq n_i,$$

and from the data pairs $(X_{ij}, Y_{ij}^{**})$, compute the double-bootstrap version $\widehat{\Theta}_i^{**\text{BLUP}}$ of $\widehat{\Theta}_i^{\text{BLUP}}$. Define

$$\widehat{\text{MSE}}_i^* = E\{(\widehat{\Theta}_i^{**\text{BLUP}} - \Theta_i^{**})^2 \mid \mathcal{X}, \mathcal{Z}^*\},$$

where $\Theta_i^{**} = \hat{\mu}^* + \underline{X}_i'\hat{\beta}^* + U_i^{**}$. Then $\widehat{\text{MSE}}_i^*$ is the direct bootstrap analogue of $\widehat{\text{MSE}}_i$. The bias of $\widehat{\text{MSE}}_i$ is estimated by

$$(2.24) \qquad \widehat{\text{bias}}_i = E(\widehat{\text{MSE}}_i^* \mid \mathcal{Z}) - \widehat{\text{MSE}}_i,$$

and a simple bias-corrected estimator is

$$(2.25) \qquad \widehat{\mathrm{MSE}}_i^{\mathrm{bc}} = \widehat{\mathrm{MSE}}_i - \widehat{\mathrm{bias}}_i = 2\widehat{\mathrm{MSE}}_i - E(\widehat{\mathrm{MSE}}_i^* \mid \mathcal{Z}).$$

See Section 3 for discussion of other approaches.

The bootstrap estimators $\widehat{\mathrm{bias}}_i$ and $\widehat{\mathrm{MSE}}_i^{\mathrm{bc}}$ are analogues of the analytical estimators $\widetilde{\mathrm{bias}}_i$ and $\widetilde{\mathrm{MSE}}_i^{\mathrm{bc}}$, respectively, introduced in Section 2.4. We shall show in Section 4.2 that the bootstrap estimators have the same orders of accuracy as their analytical counterparts, in that

$$(2.26) \qquad \widehat{\mathrm{bias}}_i = \widetilde{\mathrm{bias}}_i + O_p(n^{-2}), \qquad E(\widehat{\mathrm{bias}}_i) = E(\widetilde{\mathrm{bias}}_i) + O(n^{-2}),$$

$$(2.27) \quad \widehat{\mathrm{MSE}}_i^{\mathrm{bc}} = \widetilde{\mathrm{MSE}}_i^{\mathrm{bc}} + O_p(n^{-2}), \qquad E(\widehat{\mathrm{MSE}}_i^{\mathrm{bc}}) = E(\widetilde{\mathrm{MSE}}_i^{\mathrm{bc}}) + O(n^{-2}).$$

2.6. *Distributions* $D(z_2, z_4)$. The simplest example of a distribution $D(1, p^{-1})$ of a random variable $Z$ is perhaps the three-point distribution,

$$(2.28) \qquad P(Z = 0) = 1 - p, \qquad P(Z = \pm p^{-1/2}) = \tfrac{1}{2}p,$$

where $0 < p < 1$. Here, $E(Z) = 0$, $E(Z^2) = 1$ and $E(Z^4) = p^{-1}$. Therefore we may take $D(z_2, z_4)$ to be the distribution of $z_2^{1/2} Z$ when $p = z_2^2 / z_4$.

The Pearson family of distributions has the potential for fitting the first four moments. If (a) the first and third moments are zero, (b) the second is $z_2 = 1$ and (c) the fourth is $z_4 > 3$, implying that tails are heavier than those of the normal distribution, then the Pearson family distribution is rescaled Student's $t$. The number of degrees of freedom, $r$, is not necessarily an integer, and is given by $z_4 = 3(r-2)/(r-4)$.

Section 3 reports results of a simulation study where both the three-point and Student's $t$ distributions are used. While Student's $t$ can be employed only when kurtosis is positive, this is the case in many practical situations.

2.7. *Estimating fourth moments of $U$ and $V$*. A variety of methods can be used; the one suggested here is based on estimating moments of residuals. Define

$$W_{ij_1j_2}(s, t) = s(U_i + s_{ij_1} V_{ij_1}) + t(U_i + s_{ij_2} V_{ij_2}),$$

to which an empirical approximation is

$$\widehat{W}_{ij_1j_2}(s, t) = s(Y_{ij_1} - \hat{\mu} - X'_{ij_1}\hat{\beta}) + t(Y_{ij_2} - \hat{\mu} - X'_{ij_2}\hat{\beta}).$$

The average value, $\bar{W}_k(s, t)$, of $\widehat{W}_{ij_2}(s, t)^k$, over pairs $(j_1, j_2)$ of distinct integers $1 \le j_1, j_2 \le n_i$ and over $1 \le i \le n$, is a root-$n$ consistent estimator of the analogous average value, $w_k(s, t)$ say, of $E\{W_{ij_1j_2}(s, t)^k\}$. Now,

$$w_4(1, -1) = 2a_4 E(V^4) + 6 \frac{\sum_i (\sum_j s_{ij}^2)^2 - \sum_i \sum_j s_{ij}^4}{\sum_i n_i(n_i - 1)} (EV^2)^2,$$

where $a_4 = N^{-1} \sum_i \sum_j s_{ij}^4$. This suggests the estimator

$$\hat{\gamma}_V = \max\left[ (2a_4)^{-1} \left\{ \bar{W}_4(1,-1) - 6 \frac{\sum_i (\sum_j s_{ij}^2)^2 - \sum_i \sum_j s_{ij}^4}{\sum_i n_i(n_i-1)} \hat{\sigma}_V^4 \right\}, \hat{\sigma}_V^4 \right]$$

of $\gamma_V = E(V^4)$, which leads in turn to an estimator of $\gamma_U = E(U^4)$,

$$\hat{\gamma}_U = \max\left[ N^{-1}\left\{ \sum_{i=1}^{n} \sum_{j=1}^{n_i} (Y_{ij} - \hat{\mu} - X_{ij}'\hat{\beta})^4 \right. \right.$$

$$\left. \left. - 6\hat{\sigma}_U^2 \hat{\sigma}_V^2 \sum_{i=1}^{n} \sum_{j=1}^{n_i} s_{ij}^2 - \hat{\gamma}_V \sum_{i=1}^{n} \sum_{j=1}^{n_i} s_{ij}^4 \right\}, \hat{\sigma}_U^4 \right].$$

**3. Numerical properties.** Recall, from (2.24), that the bias of $\hat{u} = \widehat{\text{MSE}}_i$ is estimated by $\widehat{\text{bias}}_i = \hat{v} - \hat{u}$, where $\hat{v} = E(\widehat{\text{MSE}}_i^* \mid \mathcal{Z})$. The bias of $\hat{u}$ can be corrected in a broad variety of ways. Perhaps the simplest is to take $\hat{u} - \widehat{\text{bias}}_i = 2\hat{u} - \hat{v}$ as our estimator of $\text{MSE}_i$. To avoid difficulties with the sign of this quantity we might instead take as our bias correction $\hat{u} + n^{-1}g\{n(\hat{u} - \hat{v})\}$, where $g(t)$ is a smooth, symmetric, bounded function which equals $t$, or approximately $t$, when $t$ is not very far 0. One approach which incorporates this idea is to use

(3.1) $$\overline{\text{MSE}}_i = \begin{cases} \hat{u} + n^{-1}g\{n(\hat{u} - \hat{v})\}, & \text{if } \hat{u} \geq \hat{v}, \\ \hat{u}^2/[\hat{u} + n^{-1}g\{n(\hat{v} - \hat{u})\}], & \text{if } \hat{u} < \hat{v}. \end{cases}$$

The two-case definition of $\overline{\text{MSE}}_i$ ensures that this estimator is positive. The fact that the high-order correction, $\hat{u} - \hat{v}$, is captured inside the bounded function $g$ limits the detrimental effects that stochastic variation of the correction can have on overall variability, so removing the first drawback.

An elementary choice, which gives very good results in practice, is $g(t) = \text{sgn}(t)\min(|t|, nc)$, where $c$ is a positive constant. Perhaps surprisingly, $g(t) = \arctan t$ also performs well. For the sake of brevity we shall report results only for the latter estimator, although similar performance is obtained using other approaches.

In the remainder of this section we report results of a simulation study under the regression model (2.1). We took $r = 1$, $\mu = 0$, $\beta = 1$, each $n_i = 3$, $s_{ij} = 1$ for all $i$ and $j$, and $n = 60$ or 100; and we generated the $X_{ij}$'s from the Uniform distribution on $[\frac{1}{2}, 1]$. The objective was to estimate the mixed effects $\mu + \underline{X}_i\beta + U_i$. In problems of small-area estimation, this quantity can be treated as the small-area mean.

Eight different models for the distributions of $U$ and $V$ were considered, in each case centered so that both distributions had zero mean. Variances were standardized so that the ratio $\sigma_U^2/\sigma_V^2$ equaled $\frac{1}{2}$, 1 or 2, $\max(\sigma_U^2, \sigma_V^2) = 1$, and $\min(\sigma_U^2, \sigma_V^2) = \frac{1}{2}$ or 1. The models were $M_1$: $U$ and $V$ are both normal;

$M_2$: $U$ and $V$ are both $\sqrt{\chi_5^2}$; $M_3$: $U$ and $V$ are both $\chi_5^2$; $M_4$: $U$ and $V$ are both $\chi_{10}^2$; $M_5$: $U$ and $V$ are both exponential; $M_6$: $U$ is $\chi_5^2$ and $V$ is $-\chi_5^2$; $M_7$: $U$ and $V$ are both Student's $t_6$; $M_8$: $U$ and $V$ both have logistic distributions.

We used empirical measures of relative bias and coefficient of variation to quantify the performances of our methods for different distributions. Relative bias of the mean-squared error estimator was defined to be the average, over $i$, of

$$(3.2) \qquad \mathrm{RB}_i = \frac{E(\widehat{\mathrm{MSE}}_i) - \mathrm{SMSE}_i}{\mathrm{SMSE}_i},$$

$i = 1, \ldots, n$, where $E(\widehat{\mathrm{MSE}}_i)$ was estimated empirically as the average of values of $\widehat{\mathrm{MSE}}_i$ over replicates. (We shall also report the average of the absolute values of the $\mathrm{RB}_i$'s.) Likewise, $\mathrm{SMSE}_i$ was defined as the average value of $(\widehat{\Theta}_i - \Theta_i)^2$ over replicates. The coefficient of variation of the MSE estimator was taken to be the average, over $i$, of

$$\mathrm{CV}_i = \frac{\{E(\widehat{\mathrm{MSE}}_i - \mathrm{SMSE}_i)^2\}^{1/2}}{\mathrm{SMSE}_i},$$

$i = 1, \ldots, n$, where $E(\widehat{\mathrm{MSE}}_i - \mathrm{SMSE}_i)^2$ was computed by averaging $E(\widehat{\mathrm{MSE}}_i - \mathrm{SMSE}_i)^2$ over replicates.

Table 1 reports results in the case $\sigma_U^2/\sigma_V^2 = 1$, and Table 2 gives results for $\sigma_U^2/\sigma_V^2 = \frac{1}{2}$ and 2. For comparison, results for the "naive" mean-squared prediction error estimator, without any bias correction, are reported in the column headed RBN. The naive estimator is obtained by replacing $\sigma_U^2$ and $\sigma_V^2$, in the formula at (2.12), by $\hat{\sigma}_U^2$ and $\hat{\sigma}_V^2$, respectively.

In the problem of estimating predictive mean-squared error, the naive estimator is notoriously optimistic; it is significantly negatively biased. Erring by giving a falsely positive impression of reliability can significantly affect the level of debate about policy decisions based on predictions. Ideally, bias correction should remove much of this effect, producing estimators that tend to err on the side of overestimation of variance, and, against that background, to reduce the overall magnitude of bias. The results in Tables 1 and 2 show that, to a substantial extent, our bias-corrected estimator achieves this goal.

In Table 1, the average relative bias, across all models, is small, less than 10% and, in some individual cases, less than 5%. The three-point distribution tends to give lower relative bias than Student's $t$ in the case of skewed distributions, although it has slightly higher coefficient of variation. For our method, both the relative bias and the coefficient of variation tend to decrease as $n$ increases.

In marked contrast, the naive estimator of mean-squared error suffers from substantial underestimation, in the range 8%–20%. Indeed, the percentages of cases where underestimation occurred, for models $M_1$ to $M_8$, respectively, and for the pair (bias-corrected estimator, naive estimator), are $(3, 26)$, $(5, 56)$, $(26, 67)$, $(20, 63)$, $(13, 59)$, $(9, 29)$, $(3, 52)$ and $(0, 73)$, respectively. The average percentages of absolute bias, measured in terms of (median, mean), are $(12.6, 15.9)$ for our bias-corrected estimator, and $(18.8, 24.4)$ for the naive method. All these results are for the case of moment-matching using the three-point distribution.

In Table 2, to save space we give results only for models $M_3$ and $M_7$. It can be seen from those results that, for unequal variance components, both the relative bias and the coefficient of variation tend to take higher values, compared to the equal-variance case. However, the difference is not large. The three-point distribution, used to match moments, tends to give slightly better results here than the Student's $t$ approach.

The performance of normal-theory bias corrections applied to nonnormal data is well documented. For example, in the case of exponential data, use of normal theory can result in relative bias of 19% [27]. The extent of overestimation evidenced in Tables 1 and 2 is common in work on bias-correction in related problems; see the simulation results of [24, 27, 32].

TABLE 1
*Relative bias and empirical coefficient of variation under different models*

| Model | $n = 60$ | | | | | $n = 100$ | | | | |
|---|---|---|---|---|---|---|---|---|---|---|
| | 3pt | | $t$ | | | 3pt | | $t$ | | |
| | RB | CV | RB | CV | RBN | RB | CV | RB | CV | RBN |
| $M_1$ | 0.088 | 0.250 | 0.084 | 0.244 | −0.147 | 0.082 | 0.238 | 0.078 | 0.142 | −0.150 |
| | 0.091 | 0.290 | 0.100 | 0.286 | −0.131 | 0.098 | 0.280 | 0.080 | 0.162 | −0.142 |
| $M_2$ | 0.062 | 0.262 | 0.099 | 0.253 | −0.185 | 0.058 | 0.247 | 0.081 | 0.248 | −0.181 |
| | 0.089 | 0.289 | 0.103 | 0.291 | −0.187 | 0.092 | 0.286 | 0.088 | 0.274 | −0.142 |
| $M_3$ | 0.066 | 0.292 | 0.097 | 0.271 | −0.200 | 0.040 | 0.262 | 0.053 | 0.264 | −0.191 |
| | 0.095 | 0.331 | 0.101 | 0.323 | −0.200 | 0.067 | 0.298 | 0.048 | 0.301 | −0.191 |
| $M_4$ | 0.064 | 0.272 | 0.062 | 0.258 | −0.125 | 0.039 | 0.254 | 0.064 | 0.258 | −0.103 |
| | 0.076 | 0.312 | 0.099 | 0.305 | −0.121 | 0.051 | 0.279 | 0.076 | 0.289 | −0.117 |
| $M_5$ | 0.088 | 0.360 | 0.103 | 0.331 | −0.141 | 0.070 | 0.295 | 0.090 | 0.278 | −0.142 |
| | 0.108 | 0.375 | 0.111 | 0.375 | −0.163 | 0.079 | 0.327 | 0.100 | 0.315 | −0.158 |
| $M_6$ | 0.006 | 0.283 | 0.109 | 0.282 | −0.125 | 0.044 | 0.276 | 0.080 | 0.281 | −0.112 |
| | 0.075 | 0.317 | 0.121 | 0.316 | −0.125 | 0.064 | 0.312 | 0.088 | 0.313 | −0.112 |
| $M_7$ | 0.100 | 0.331 | 0.099 | 0.287 | −0.158 | 0.028 | 0.262 | 0.015 | 0.246 | −0.100 |
| | 0.106 | 0.376 | 0.099 | 0.327 | −0.166 | 0.036 | 0.280 | 0.036 | 0.268 | −0.115 |
| $M_8$ | 0.104 | 0.299 | 0.065 | 0.260 | −0.112 | 0.093 | 0.281 | 0.056 | 0.244 | −0.080 |
| | 0.100 | 0.326 | 0.119 | 0.318 | −0.140 | 0.097 | 0.288 | 0.066 | 0.277 | −0.114 |

The first line in each row gives median values, and the second line, means.

TABLE 2

*Relative bias and empirical coefficient of variation for models* $M_3$
*and* $M_7$

|  | Model | **3pt** | | **$t$** | | |
|---|---|---|---|---|---|---|
|  |  | **RB** | **CV** | **RB** | **CV** | **RBN** |
| $\sigma_U^2/\sigma_V^2 = \frac{1}{2}$ | $M_3$ | 0.110 | 0.099 | 0.308 | 0.309 | −0.181 |
|  |  | 0.103 | 0.081 | 0.324 | 0.346 | −0.190 |
|  | $M_7$ | 0.124 | 0.100 | 0.289 | 0.304 | −0.135 |
|  |  | 0.109 | 0.114 | 0.338 | 0.344 | −0.135 |
| $\sigma_U^2/\sigma_V^2 = 2$ | $M_3$ | 0.099 | 0.104 | 0.310 | 0.307 | −0.118 |
|  |  | 0.112 | 0.111 | 0.358 | 0.358 | −0.140 |
|  | $M_7$ | 0.105 | 0.081 | 0.326 | 0.323 | −0.116 |
|  |  | 0.111 | 0.099 | 0.379 | 0.366 | −0.137 |

The first line in each row gives median values, and the second line,
means.

Finally we report on a comparison of our method with the parametric
jackknife approach suggested by Jiang, Lahiri and Wan [19]. The latter is
awkward to implement unless there is a closed-form expression for the lead-
ing term in an expansion of mean-squared prediction error, as a function of
unknown parameters. Among the models $M_1$–$M_8$, a closed-form expression
exists only for the first (i.e., normal–normal) model. Moreover, only in this
case does the best predictor (the small-area estimator in the case of the
jackknife) have a closed-form expression.

In the normal–normal model, and when $\sigma_U^2/\sigma_V^2 = 1$ and $n = 60$, the me-
dian (mean) relative bias and the median (mean) coefficient of variation are
0.035 (0.049) and 0.262 (0.298), respectively. When $n = 100$ the correspond-
ing values are 0.034 (0.047) and 0.156 (0.182). For unequal variance ratios
the relative biases are close to these values, while the coefficients of variation
are higher.

Comparing these results with those in Table 1, it can be seen that the
jackknife method, which uses full knowledge of the error distributions, per-
forms better in terms of relative bias but is inferior in terms of coefficient of
variation, relative to the nonparametric bootstrap method. The impact of
deviation from normality of error distributions has been reported by Prasad
and Rao [27] and Wang and Fuller [32].

## 4. Theoretical properties.

4.1. *Rigorous formulations of* (2.11), (2.14), (2.15), (2.19) *and* (2.21).
We begin by stating, and discussing, regularity conditions. Of the $n_i$'s, $s_{ij}$'s,

$X_{ij}$'s and distributions of $U$ and $V$ we ask that:

(4.1)

(a) $\sup_i n_i < \infty$ and each $n_i \geq 2$, (b) $C_1 \leq s_{ij} \leq C_2$ for constants $0 < C_1 < C_2 < \infty$ and for all $i$ and $j$, (c) the vectors $X_{ij}$ are conditioned-upon values of independent copies of the random $r$-vector $X$, the distribution of which is continuous and satisfies $P(\|X\| \leq C_3) = 1$ for some $0 < C_3 < \infty$, (d) all moments of $U$ and $V$ are finite, $E(U) = E(V) = 0$ and $\sigma_V > 0$, and (e) the eigenvalues of the $n_i \times n_i$ matrix $\mathbf{T}_i$ are bounded away from zero, uniformly in $1 \leq i \leq n < \infty$.

The conditions on $n_i$ in (a) do not require any constraints on the long-run frequencies of different values of the $n_i$'s. As a result, the functions $\psi_j$ in the expansions in Section 2.4, which depend on $n$, may not converge as $n \to \infty$. However, they will converge if we assume in addition that the proportion of values of $i$, $1 \leq i \leq n$, for which $n_i$ takes any given value between 2 and $\sup_i n_i$, converges as $n \to \infty$. Nevertheless, without this condition the functions $\psi_j$ are uniformly bounded.

Condition (c), in (4.1) and on the $X_{ij}$'s, can be weakened, and in particular it is not essential to assume that the distribution of each $X_{ij}$ is the same for all $i$ and $j$. However, without that constraint, more complex assumptions have to be made in order to ensure that the distributions of the $X_{ij}$'s do not become "asymptotically degenerate" as $n \to \infty$. If this occurs, then it could adversely affect assumptions made in Section 2.3 about the rank of the matrices $\mathbf{P}$ and $\bar{\mathbf{P}}$; those assumptions automatically hold, with probability 1 with respect to the process generating the $X_{ij}$'s, under the present conditions. Concerning assumption (d), it is not essential to assume that $\sigma_U > 0$. Assumption (e) is a restriction on choice of the $s_{ij}$'s.

As noted at the end of Section 2.4, in general it is necessary to introduce a ridge parameter to ensure that $\widehat{W}_i$ is nonsingular. We do this by replacing $\mathrm{SSE}_1$ by $\max(\mathrm{SSE}_1, B_1 n^{-B_2})$, for some $B_1 > 0$ and $B_2 \geq 2$, in the definition of $\mathrm{SSE}_1$. It will be assumed below that this has been done. Depending on the distributions of $U$ and $V$, this can slightly alter the definitions of $\hat{\sigma}_V^2$ and $\hat{\sigma}_U^2$. The ridge parameter is not necessary for Theorem 1 if the distribution of $V$ is absolutely continuous.

Recall from Section 2.4 that $\xi_0 = (\sigma_U^2, \sigma_V^2)$, $\xi_1 = (\sigma_U^2, \sigma_V^2, \gamma_U, \gamma_V)$ and $\psi_0(\xi_0) = \sigma_U^2 a_i^{-1} \sigma_V^2/(\sigma_U^2 + a_i^{-1}\sigma_V^2)$, where $\gamma_U = E(U^4)$ and $\gamma_V = E(V^4)$. Here, and in (4.2)–(4.6) below, we suppress the dependence of the functions $\psi_0$, $\psi_1$ and $\psi_2$ on $i$, and expectations are interpreted as conditional on $\mathcal{X}$.

THEOREM 1. *If (4.1) holds, then, for a class of realizations of $\mathcal{X}$ that arises with probability 1, and for $k = 1, 2$,*

$$\mathrm{MSE}_i = \frac{\sigma_U^2 a_i^{-1} \sigma_V^2}{\sigma_U^2 + a_i^{-1}\sigma_V^2} + n^{-1}\psi_1(\xi_1) + O(n^{-2}),$$
(4.2)

(4.3)        $E\{\psi_0(\hat{\sigma}_U^2, \hat{\sigma}_V^2)\} = \psi_0(\xi_0) + n^{-1}\psi_2(\xi_1) + O(n^{-2}),$

(4.4)    $E\{\psi_k(\hat{\sigma}_U^2, \hat{\sigma}_V^2, \hat{\gamma}_V, \hat{\gamma}_V)\} = \psi_k(\xi_1) + O(n^{-1}),$

*uniformly in $1 \le i \le n$, where the functions $\psi_1$ and $\psi_2$ are determined solely by the design variables $X_{ij}$ and weights $s_{ij}$ for $1 \le i \le j \le n$, depend on $n$, are bounded in a neighborhood of $\xi_1$, and are infinitely differentiable.*

A proof of (4.2) is given in the web version of this paper [17]. Derivations of (4.3) and (4.4) are similar but simpler. Together, (4.2)–(4.4) imply (2.21), which asserts that $\widehat{\mathrm{MSE}}_i^{\mathrm{bc}}$, defined at (2.20), has bias equal to $O(n^{-2})$.

4.2. *Theory for the bootstrap.* In Sections 2.4 and 2.5 we discussed analytical and bootstrap-based bias corrections, respectively. In particular, $\widehat{\mathrm{MSE}}_i$, at (2.16), was an analytical estimator of $\mathrm{MSE}_i$; the associated analytical bias estimator was $\widehat{\mathrm{bias}}_i$, at (2.22); and the bias-corrected estimator was $\widehat{\mathrm{MSE}}_i^{\mathrm{bc}} = \widehat{\mathrm{MSE}}_i - \widehat{\mathrm{bias}}_i$, at (2.20). In the same vein, $\widetilde{\mathrm{MSE}}_i$, at (2.23), was a bootstrap estimator of $\mathrm{MSE}_i$; the corresponding bootstrap estimator of bias was $\widetilde{\mathrm{bias}}_i$, at (2.24); and the resulting bias-corrected estimator was $\widetilde{\mathrm{MSE}}_i^{\mathrm{bc}}$, at (2.25).

The effectiveness of the bootstrap approach is reflected in the fact that it gives a degree of correction that is identical to that provided by the analytical method, up to terms of order $n^{-2}$, as the next result shows. For definiteness we assume there that the moment-matching bootstrap method is based on the three-point distribution at (2.28); the Student's $t$ model does not permit correction for negative kurtosis. We suppose too that the ridge parameter defined two paragraphs above Theorem 1 is incorporated into the definition of $\hat{\sigma}_V^2$. (This turned out not to be necessary in our simulation study, even though the three-point distribution has positive mass at zero. That can be explained by noting that the probability of difficulty being caused by the positive mass at zero is exponentially small, as a function of sample size, whereas we used only polynomially many bootstrap simulations.)

THEOREM 2. *If (4.1) holds, and if the distribution at (2.28) is used to implement the moment-matching bootstrap, then for a class of realizations of $\mathcal{X}$ that arises with probability 1,*

(4.5)    $\widehat{\mathrm{MSE}}_i - \widetilde{\mathrm{MSE}}_i = O_p(n^{-2}), \qquad E\{\widehat{\mathrm{MSE}}_i - \widetilde{\mathrm{MSE}}_i\} = O(n^{-2}),$

(4.6)    $\widehat{\mathrm{bias}}_i - \widetilde{\mathrm{bias}}_i = O_p(n^{-2}), \qquad E\{\widehat{\mathrm{bias}}_i - \widetilde{\mathrm{bias}}_i\} = O_p(n^{-2}),$

*uniformly in $1 \le i \le n$.*

Results (4.2)–(4.4) established the efficacy of the analytical approach to bias correction. In combination with those properties, (4.5) and (4.6) do the same for the bootstrap approach, by establishing (2.26) and (2.27). A proof of Theorem 2 is given by Hall and Maiti [17].

## REFERENCES

[1] BATTESE, G. E., HARTER, R. M. and FULLER, W. A. (1988). An error-components model for prediction of county crop areas using survey and satellite data. *J. Amer. Statist. Assoc.* **83** 28–36.

[2] BELL, W. (2001). Discussion of "Jackknife in the Fay–Herriot model with an application," by Jiang, Lahiri, Wan and Wu. In *Proc. Seminar on Funding Opportunity in Survey Research* 98–104. Council of Professional Associations on Federal Statistics, Washington.

[3] BOOTH, J. G. and HOBERT, J. P. (1998). Standard errors of prediction in generalized linear mixed models. *J. Amer. Statist. Assoc.* **93** 262–272. MR1614632

[4] CARROLL, R. J. and HALL, P. (1988). Optimal rates of convergence for deconvolving a density. *J. Amer. Statist. Assoc.* **83** 1184–1186. MR0997599

[5] CHEN, S. and LAHIRI, P. (2003). A comparison of different MSPE estimators of EBLUP for the Fay–Herriot model. *Proc. Survey Research Methods Section* 905–911. Amer. Statist. Assoc., Alexandria, VA.

[6] DAS, K., JIANG, J. and RAO, J. N. K. (2004). Mean squared error of empirical predictor. *Ann. Statist.* **32** 818–840. MR2060179

[7] DATTA, G. S. and GHOSH, M. (1991). Bayesian prediction in linear models: Applications to small area estimation. *Ann. Statist.* **19** 1748–1770. MR1135147

[8] DATTA, G. S. and LAHIRI, P. (2000). A unified measure of uncertainty of estimated best linear unbiased predictors in small area estimation problems. *Statist. Sinica* **10** 613–627. MR1769758

[9] DELAIGLE, A. and GIJBELS, I. (2004). Practical bandwidth selection in deconvolution kernel density estimation. *Comput. Statist. Data Anal.* **45** 249–267. MR2045631

[10] DOMÍNGUEZ, M. A. and LOBATO, I. N. (2003). Testing the martingale difference hypothesis. *Econometric Rev.* **22** 351–377. MR2018885

[11] EL-AMRAOUI, A. and GOFFINET, B. (1991). Estimation of the density of $G$ given observations of $Y = G + E$. *Biometrical J.* **33** 347–355. MR1129072

[12] FAN, J. (1991). On the optimal rates of convergence for nonparametric deconvolution problems. *Ann. Statist.* **19** 1257–1272. MR1126324

[13] FAN, J. (1992). Deconvolution with supersmooth distributions. *Canad. J. Statist.* **20** 155–169. MR1183078

[14] FAN, Y. and LI, Q. (2002). A consistent model specification test based on the kernel sum of squares of residuals. *Econometric Rev.* **21** 337–352. MR1944979

[15] FLACHAIRE, E. (2002). Bootstrapping heteroskedasticity consistent covariance matrix estimator. *Comput. Statist.* **17** 501–506. MR1952694

[16] GONZÁLEZ MANTEIGA, W., MARTÍNEZ MIRANDA, M. D. and PÉREZ GONZÁLEZ, A. (2004). The choice of smoothing parameter in nonparametric regression through wild bootstrap. *Comput. Statist. Data Anal.* **47** 487–515. MR2100562

[17] HALL, P. and MAITI, T. (2005). Nonparametric estimation of mean-squared prediction error in nested-error regression models. Available at http://arxiv.org/abs/math/0509493.

[18] HARVILLE, D. A. and JESKE, D. R. (1992). Mean squared error of estimation or prediction under a general linear model. *J. Amer. Statist. Assoc.* **87** 724–731. MR1185194

[19] Jiang, J., Lahiri, P. and Wan, S.-M. (2002). A unified jackknife theory for empiri-
cal best prediction with $M$-estimation. *Ann. Statist.* **30** 1782–1810. MR1969450

[20] Kacker, R. and Harville, D. A. (1984). Approximations for standard errors of
estimators of fixed and random effects in mixed linear models. *J. Amer. Statist.
Assoc.* **79** 853–862. MR0770278

[21] Kauermann, G. and Opsomer, J. D. (2003). Local likelihood estimation in gener-
alized additive models. *Scand. J. Statist.* **30** 317–337. MR1983128

[22] Lahiri, P. (2003). On the impact of bootstrap in survey sampling and small-area
estimation. *Statist. Sci.* **18** 199–210. MR2019788

[23] Lahiri, P. (2003). A review of empirical best linear unbiased prediction for the Fay–
Herriot small-area model. *Philippine Statistician* **52** 1–15.

[24] Lahiri, P. and Rao, J. N. K. (1995). Robust estimation of mean squared error of
small area estimators. *J. Amer. Statist. Assoc.* **90** 758–766. MR1340527

[25] Li, Q., Hsiao, C. and Zinn, J. (2003). Consistent specification tests for semiparamet-
ric/nonparametric models based on series estimation methods. *J. Econometrics*
**112** 295–325. MR1951146

[26] Li, T. and Vuong, Q. (1998). Nonparametric estimation of the measurement error
model using multiple indicators. *J. Multivariate Anal.* **65** 139–165. MR1625869

[27] Prasad, N. G. N. and Rao, J. N. K. (1990). The estimation of mean squared error
of small-area estimators. *J. Amer. Statist. Assoc.* **85** 163–171. MR1137362

[28] Prášková, Z. (2003). Wild bootstrap in RCA(1) model. *Kybernetika* (*Prague*) **39**
1–12. MR1980120

[29] Rao, J. N. K. (2003). *Small Area Estimation.* Wiley, Hoboken, NJ. MR1953089

[30] Rao, J. N. K. and Choudhry, G. H. (1995). Small area estimation: Overview and
empirical study. In *Business Survey Methods* (B. G. Cox, D. A. Binder, B. N.
Chinnappa, A. Christianson, M. J. Colledge and P. S. Kott, eds.) 527–542. Wiley,
New York.

[31] Stukel, D. M. and Rao, J. N. K. (1997). Estimation of regression models with
nested error structure and unequal error variances under two and three stage
cluster sampling. *Statist. Probab. Lett.* **35** 401–407. MR1483027

[32] Wang, J. and Fuller, W. A. (2003). The mean squared error of small area pre-
dictors constructed with estimated area variances. *J. Amer. Statist. Assoc.* **98**
716–723. MR2011685

Centre for Mathematics                          Department of Statistics
  and Its Applications                          Iowa State University
Australian National University                  221 Snedecor Hall
Canberra, ACT 0200                              Ames, Iowa 50011
Australia                                       USA
E-mail: Peter.Hall@anu.edu.au                   E-mail: taps@iastate.edu